\documentclass[11pt]{article}
\usepackage{algorithm}
\usepackage{algorithmic}
\usepackage[margin=0.8in]{geometry}
\usepackage{pgfplots}
\usepackage{amsmath}
\usepackage{amssymb}
\usepackage{algorithmic}
\usepackage{algorithm}
\usepackage{mathtools}
\usepackage{color}
\usepackage{url}
\usepackage{tikz}
\usetikzlibrary{arrows,
	arrows.meta,
	bending,calc,patterns,angles,quotes}
\usepackage{float}

\DeclareMathOperator*{\argmax}{argmax}

\newcommand{\var}{\mbox{Var}}

\makeatletter
\providecommand*{\cupdot}{%
	\mathbin{%
		\mathpalette\@cupdot{}%
	}%
}
\newcommand*{\@cupdot}[2]{%
	\ooalign{%
		$\m@th#1\cup$\cr
		\hidewidth$\m@th#1\cdot$\hidewidth
	}%
}
\makeatother

\newtheorem{example}{Example}
\newtheorem{conjecture}{Conjecture}

\newtheorem{theorem}{Theorem}

\newtheorem{remark}{Remark}
\newtheorem{definition}{Definition}

\def\endpf{{\ \hfill\hbox{\vrule width1.0ex height1.0ex}\parfillskip 0pt
	}}
	\newenvironment{proof}{\noindent{\bf Proof:}}{\endpf}
\begin{document}
	\title{Minimizing the externalities variance in a LCFS-PR $M/G/1$  queue under various constraints}
	\author{Royi Jacobovic\footnote{The Department of Statistics and Operations Research, School of Mathematical Sciences, Tel Aviv University, Israel; {\tt royijacobo@tauex.tau.ac.il}}$\ $  and Nikki Levering\thanks{Korteweg-de Vries Institute; University of Amsterdam; 1098 XG Amsterdam; Netherlands. 
			{\tt n.a.c.levering@uva.nl}} }
	\date{\today}
	\maketitle
	\begin{abstract} 
Consider a LCFS-PR $M/G/1$ queue and assume that at time $t = 0$, there are $n+2$ customers $c_1,c_2,...,c_{n+1},c$ who arrived in that order such that $t = 0$ is the arrival time of $c$. Then, the externalities which are generated by $c$ is the total waiting time that would be saved by $c_1,c_2,...,c_{n+1}$ if $c$ reduced his service requirement to zero. Motivated by some applications, this work is about the minimization of the externalities variance under various constraints. \end{abstract}
	
	
    \bigskip
	\noindent {\bf AMS Subject Classification (MSC2010):}  60K25, 60K30, 90C08, 90C25, 90C27.

\section{Introduction}\label{sec: introduction}

\subsection{Externalities in queues}
Congestion is a phenomenon which characterizes storage systems such as
telecommunication networks, production-lines, road-traffic, etc. In practice,
excessive congestion leads to reduced system efficiency. Motivated by
this issue, queueing theory is the study of  mathematical modeling and
optimization of congested systems which are operated in various stochastic
environments. 

An observation which dates back to classical works, e.g., \cite{Edelson1975,Leeman1964,Naor1969}, is as follows: Each customer who joins the queue will cause a potential delay to the customers who will be waiting during his service periods. Moreover, the individual utility maximization of any joining customer does not take into account the misery which he will cause to those who will have to wait because of him. Therefore, a plausible remedy to this market-failure is to impose a toll on customers who are joining the queue.  

Basically, the social-planner problem is how to determine this toll in order to maximize the long-run average total welfare of the customers (for simplicity, assume that there are no costs due to service provision). One natural pricing policy to think about is called \textit{internalization of externalities}, i.e., make each joining customer pay the expected loss which he will cause to others during his sojourn in the system, when all other customers are getting their service according to the socially optimal resource allocation. Some studies about this idea and its variations to queues in which the customers have no influence on their service times are, e.g., \cite{Carlin1970,Dewan1990,Hassin1985,Hassin1994, Mendelson1990}. In addition, there are also some works about this idea in the context of queues in which customers have influence on their service times, e.g., \cite{Ha2001,Haviv2014,Jacobovic2022a,Jacobovic2022b}. Another genre of problems which are related to the notion of externalities concerns server-allocation (scheduling) problems in multi-class queues (especially with dynamic class types) \cite{Chan2021,Hu2022,Huang2015,Liu2022}. For a more detailed discussion about the relation between these works and the notion of externalities, see  \cite[Section 1.1]{Jacobovic2022}.  
 
Clearly, the externalities which are generated by a tagged customer, i.e., the total waiting time  which would be saved by the other customers if the tagged customer reduced his service requirement to zero and left the system, is a random variable. Notably, since this random variable is very applicable, it makes sense to study its distribution. The prominent works in this direction are \cite{Haviv1998,Jacobovic2022,Jacobovic2023}, and, as to be explained in the sequel, the current work adds to this strand of literature.     

\subsection{Externalities in LCFS-PR M/G/1 queue}
Consider a single-server queue with customers who arrive 
with respect to 
a Poisson process with rate $\lambda>0$. Their service demands are iid nonnegative random variables which are distributed according to $F(\cdot)$ and  independent from the arrival process. For simplicity of notation, for each $i\geq1$, let $\mu_i$ be the $i$-th moment of $F(\cdot)$ and assume that the system is stable, i.e., $\rho\equiv\lambda\mu_1<1$. Also, let the system 
be operated according to a preemptive last-come, first-served (LCFS-PR) discipline. In addition, assume that at time $t=0$, there exist $n+1\geq1$ customers $c_1,c_2,\ldots , c_{n+1}$ in the system such that for each $1\leq i\leq n+1$, the remaining service time of $c_i$ is equal to $v_i$. For clarity and without loss of generality, for every $1\leq i<j\leq n+1$ we shall assume that $c_i$ arrived before $c_j$. 

Now, assume that there is an additional customer $c$ who arrives at time $t=0$ and whose service demand is equal to $x>0$. The \textit{externalities} which are created by $c$ is the total waiting time that $c_1,c_2,\ldots,c_{n+1}$ would save if $c$ reduced his service demand to zero and left the system immediately. Importantly, the externalities are expressed via a random variable $\mathcal{E}$ with a distribution which is determined by $\lambda,n,v_1,v_2,\ldots,v_{n+1},x$ and $F(\cdot)$. 

In a recent study \cite[Eq. (19)]{Jacobovic2023} the mean and variance of $\mathcal{E}$ were derived and their formulae are as follows:
\begin{equation}\label{eq: mean}
    E\mathcal{E}=\frac{nx}{1-\rho}\,,
\end{equation}
\begin{equation}\label{eq: variance}
\var\left(\mathcal{E}\right)
=\frac{\lambda\mu_2}{(1-\rho)^3}\left[nx+2\sum_{1\leq k< \ell\leq n}\left(x-\sum_{i=k}^\ell v_i\right)^+\right]\,.
\end{equation}
It is very eminent that both of the above-mentioned formulae are independent of $v_{n+1}$. In addition, it is also a bit surprising that the expected externalities are invariant with respect to the values of $v_1,v_2,\ldots,v_{n+1}$. 

\subsection{Problem statement}
Assume that the system has a manager who knows the statistical assumptions of the model. In addition, assume that, at time $t=0$, he observes the vector 
\begin{equation}
    \left(x, v_{n+1}, x+\sum_{i=1}^{n+1}v_i, n\!+\!1\right)\,,
\end{equation} which includes four coordinates: 
\begin{enumerate}
    \item $x$ which is the remaining service time of the customer $c$ who has just entered the service position.

    \item $v_{n+1}$ which is the remaining service time of $c_{n+1}$ who has just left the service position.

    \item $x+\sum_{i=1}^{n+1}v_i$ which equals the total workload in the system at $t=0$.

    \item $n+1$ which is the number of customers who exist in the system at time $t=0$.
\end{enumerate} In particular, this information scheme looks reasonable in settings in 
which the manager has limited memory capacity, and he is unable to store the remaining
service times of all customers.

With this information, the manager may compute $E\mathcal{E}$ via \eqref{eq: mean} but he is unable to compute $\text{Var}\left(\mathcal{E}\right)$ via \eqref{eq: variance}. Thus, the manager might be willing to compute the range of possible variance values given his available information at $t=0$, i.e., he needs to derive the infimum and supremum of $\text{Var}\left(\mathcal{E}\right)$ over all parameterizations which are consistent with his available information at time $t=0$. Note that we refer to an infimum (resp. a supremum) problem and not to a minimum (resp. maximum) problem because the manager knows that there are $n+2$ customers in the system. Thus, any parameterization which is consistent with this knowledge must be such that $v_i>0$ for every $1\leq i\leq n$. In this work we focus on the infimum problem, as the supremum is evident. That is, 
since $v_1$ and $v_n$ are both contributing
to the least number of sums in \eqref{eq: variance},
it is clear that $\var\left(\mathcal{E}\right)$ reaches
its supremum if $\left(v_1,v_2,\ldots,v_n\right)$ is either $(w,0,0,\ldots,0)$ or $(0,0,\ldots,0,w)$ where $w\equiv\sum_{i=1}^nv_i$.

Eminently, when $w\geq nx$ the infimum problem is simple and hence we consider the next setup: Let $n\in\mathbb{N},\  x>0$, and $w\in(0,nx)$ be fixed parameters and define a function
\begin{equation}\label{eq: f definition}
f(\boldsymbol{v})\equiv\sum_{1\leq k\leq \ell\leq n}\left(x-\sum_{i=k}^\ell v_i\right)^+\ \ , \ \ \boldsymbol{v}=(v_1,v_2,\ldots,v_n)\in\mathbb{R}^n_+\,.
\end{equation}
Then, the infimum problem is translated to the next minimization: 
\begin{equation}\label{opt: continuous}
\min_{\boldsymbol{v}\in\Lambda}f(\boldsymbol{v})\,,
\end{equation}
where
\begin{equation}\label{eq: lambda}
    \Lambda\equiv\left\{\boldsymbol{v}\in\mathbb{R}^n_+\ ; \ \sum_{i=1}^nv_i\leq w\right\}\,.
\end{equation}

\subsection{Combinatorial version}
Now, we present a combinatorial version of the continuous minimization \eqref{opt: continuous}. Specifically, it is stated as follows: With the parameters of the original problem, denote 
$m\equiv\left\lfloor \frac{w}{x}\right\rfloor$ as well as $r\equiv w-mx$. Then, define the domain 
\begin{equation}
\Upsilon\equiv\left\{\boldsymbol{v}\in\mathbb{R}^n_+\ ;\ \sum_{i=1}^nv_i \leq w\ , \ \left|\left\{i;v_i=x\right\}\right|=m\ , \ \left|\left\{i;v_i=r\right\}\right|\geq1\right\}\,.
\end{equation} 
Verbally, if $r>0$, then $\Upsilon $ is the collection of all real $n$-dimensional vectors having $m$ components which equal $x$, one component equals $r$, and all other components equal zero. Otherwise, $w=mx$ (i.e., $r=0$), and $\Upsilon $ is the collection of all real $n$-dimensional vectors which have $m$ components which equal $x$, and all other components equal zero. Notice that $\Upsilon $ is a finite discrete set in both cases. 

The  combinatorial minimization to be studied in this paper is:
\begin{equation}\label{opt: combinatorial}
    \min_{\boldsymbol{v}\in\Upsilon }f(\boldsymbol{v})\,.
\end{equation}
Furthermore, we will explore its interplay with the continuous minimization \eqref{opt: continuous}. Also, besides being a straightforward combinatorial version of the continuous minimization \eqref{opt: continuous}, we note in passing that \eqref{opt: combinatorial} might be associated with three additional interpretations which are explained below. We believe that each of these interpretations sheds some light on a different aspect of the combinatorial minimization \eqref{opt: combinatorial}. Especially, notice that the third interpretation is also relevant for the continuous minimization \eqref{opt: continuous}. 

\subsubsection{Card collection problem}
Consider a collection of $n$ cards $\mathcal{C}\equiv\left\{C_1,C_2,\ldots,C_n\right\}$. Let $m$ be an integer which is less than $n$ and define 
\begin{equation}
    \mathcal{C}_m\equiv \left\{\{C_{i_1},C_{i_2},\ldots,C_{i_m}\};1\leq i_1<i_2<\ldots<i_m\leq n\right\}\,.
\end{equation}
Any subset of successive cards is called a \textit{series}, i.e., the set of all series is
\begin{equation}
    \mathcal{S}\equiv\left\{\left\{C_k,C_{k+1},\ldots,C_\ell\right\};1\leq k\leq\ell\leq n\right\}\,.
\end{equation}
In addition, there is a single player whose task is to choose a sub-collection $S$ of $m$ cards which participates in the maximal number of series, i.e., his needs to choose $S\in\mathcal{C}_m$ which maximizes
 \begin{equation}
     \sum_{S'\in\mathcal{S}}\textbf{1}\left(S\cap S'\neq\emptyset\right)\ \ , \ \ S \in \mathcal{C}_m\,.
 \end{equation}
Notably, this problem coincides with \eqref{opt: combinatorial}  parameterized by $x=1$, $m=w$ and $r=0$. 

Observe that when  $r\in(0,1)$, \eqref{opt: combinatorial} is describing the following stochastic version of the deterministic card collection problem: The player is required to choose $m$ cards and an additional card. The sub-collection of cards which the player possesses after making his choice is random. It will include the $m$ cards for sure, but the additional card will be included with probability $r$. Then, the purpose of the player is to choose the $m$ cards and the additional card in order to maximize the expected number of series which have nonempty intersection with his cards sub-collection.

\subsubsection{Location problem}\label{subsec: covering}
\begin{figure}
    \centering
    \includegraphics{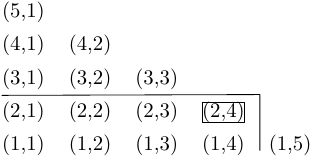}
    \caption{Example of a triangular lattice $L$ with $n = 5$.
    Station (2,4) is activated, and serves all sites in the lower left rectangle.}
    \label{fig:coveringlocation}
\end{figure}

Consider the triangular lattice 
\begin{equation}
    L\equiv\left\{(k,\ell)\in\mathbb{N}^2;\ k+\ell\leq n+1\right\}\,,
\end{equation}
and assume that each point belonging to $L$ represents a site location. In addition, assume there is a planner who should allocate $1\leq m<n$ stations on the hypotenuse of the lattice. Assume that for each $1\leq k\leq n$, the station which is located on $(k,\ell)$ provides service to all sites on
\begin{equation}
\left\{(k',\ell')\in L;k'\leq k \,, \ell' \leq \ell \right\}\,.
\end{equation}
Spatially, this means that every station provides service to all sites which are neither higher than it nor right to it (see e.g.\ Figure~\ref{fig:coveringlocation}). Then, under this assumption, the planner's objective is to pick the allocation of the stations in order to maximize the number of sites where service is available. Note that this problem is a special case of \eqref{opt: combinatorial} with the parameters $x=1$, $w=m$ and $r=0$.

For more maximal covering location problems, see a recent overview \cite{Garcia2019}. In addition, a slightly different class of covering problems is 
of the budgeted maximal coverage, see, e.g., \cite{Khuller1999}.

\subsubsection{Minimization of the excess demand for public goods}
In economic theory, a \textit{public good} is any good which has the following two properties:
\begin{enumerate}
    \item It is \textit{non-excludable}, i.e., there is no way to prevent its consumption by customers who are willing to do so. 

    \item  It is \textit{non-rivalrous}, i.e., its consumption by one customer does not limit the ability of other customers to consume it.  
\end{enumerate}
Some traditional examples of public goods are, e.g., national security, law enforcement, and clean air. In general, public goods are usually supplied by the government, since the private sector has very limited incentives to produce them (see, e.g.,  \cite{Stiglitz1977}). For more information about the theory and applications of public goods, see, e.g., \cite{Blumel1986}.

Now, consider an economy with $n$ public goods $G_1,G_2,\ldots,G_n$, and $n(n+1)/2$ customers who are doubly-indexed by $(k,\ell)$ such that $1\leq k\leq \ell\leq n$. The government has a budget of $w=mx$ monetary units and the cost of producing the allocation $\boldsymbol{v}\in\mathbb{R}^n_+$, in which for each $1\leq i\leq n$ there are $v_i$ units of the $i$-th good is $x\sum_{i=1}^nv_i$. In addition, for each $1\leq k\leq \ell\leq n$, the $(k,\ell)$-th customer demands to consume one unit of goods from the set $\{G_k,G_{k+1},\ldots,G_\ell\}$. In particular, this means that the $(k,\ell)$-th customer is indifferent between the goods $G_k,G_{k+1},\ldots,G_\ell$, i.e., he considers these goods as perfect substitutes. In order to demonstrate the market structure, one may think about $G_1,G_2,\ldots,G_n$ as food products (potentially supplied by the government) which are ordered according to their nutritional value. At the same time, $G_1,G_2,\ldots,G_n$ are reversely ordered with respect to their taste. Thus, the above-described market is such that the demand for the intermediate food products is higher than the demand for either healthy (but non-tasty) or tasty (but unhealthy) food products.   

The government's objective is to minimize the total excess-demand of the customers in the economy given its budget constraint. Formally, 
if $G_1,G_2,\ldots,G_n$ are indivisible (resp. divisible)
this problem is translated into \eqref{opt: combinatorial}  (resp. \eqref{opt: continuous}). 

Notice that, in both cases, the government's problem is how to choose the variety of products which are supplied to its citizens. For some other product variety management problems which appear in the existing literature, see, e.g.,  \cite{Lancaster1990, Ramdas2003} and the references therein.

\subsection{Additional related literature}
Technically, $f(\cdot)$ is a piecewise-linear convex function on 
the compact domain $\Lambda$, i.e., there exists an optimal solution of \eqref{opt: continuous}. Furthermore, this means that for any given parameterization of \eqref{opt: continuous} standard subgradient methods yield a sequence which converges to the numerical value of an optimal solution with adequate convergence rates. A survey about subgradient methods is given in \cite{Boyd2003}. 

In addition, observe that by using the epigraph of the convex piecewise-linear function $f(\cdot)$, the continuous minimization \eqref{opt: continuous} can be rephrased as a linear program~(LP). Therefore, it is possible to solve the continuous minimization \eqref{opt: continuous} in polynomial time, see, e.g., \cite{Karmarkar1984}. In fact, it is clear that adding the constraints $v_i\leq x$ for each $1\leq i\leq n$ will not cause any difference in the resulting minimization. As a result, \eqref{opt: continuous} may also be phrased as a special case of an \textit{interval linear programming}. Thus, the analysis in \cite{Ben-Israel1968,Ben-Israel1970,Zlobec1970,Zlobec1973} includes an explicit expression of each optimal solution of \eqref{opt: continuous} in terms of some projection matrix and the generalized inverse of the constraints matrix in the resulting linear program. 

When $r=0$, the combinatorial problem \eqref{opt: combinatorial}\label{subsec: literature combinatorial} may be considered as a version of the 0-1 knapsack problem \textit{with dependencies}. Notably, even when there are no dependencies between the items, the classic 0-1 knapsack problem is 
NP-hard \cite{Freville2004}. To the best of our knowledge, there 
is not much work about knapsack problems with dependencies. This strand of literature includes the more traditional quadratic knapsack problem.  A survey of that topic is given in \cite{Pisinger2007}. In addition, consider the recent paper \cite{Beliakov2022} about knapsack problems with dependencies through non-additive measures and Choquet integral. Other non-standard versions of the knapsack problem are summarized in  \cite{Lin1998}. 

Recall the interpretation of the combinatorial minimization \eqref{opt: combinatorial} as a card collection problem, and assume that $r=0$. Then, there is also a relation to the \textit{interval stabbing problems} which appear in \cite{Schmidt2009}. This reference
contains a discussion about the computational complexity of finding the number of series which have nonempty intersection with a given subset of $m$ cards. However, to the best of our understanding, there has been no discussion about the derivation of the optimal subset of $m$ cards, i.e., the subset which has a nonempty intersection with the maximal number of series. 

\subsection{Contribution and paper organization}
The present work includes a characterization of a set of highly-structured optimal solutions of \eqref{opt: combinatorial}. This yields an applicable numerical procedure to construct such an optimal solution \eqref{opt: combinatorial}. 
In addition, it is shown that for some families of parameterizations, e.g., when $r=0$ but not only, the continuous minimization \eqref{opt: continuous} is reduced to the combinatorial one \eqref{opt: combinatorial}. Finally, we put forth a conjecture that this reduction is true for every parameterization.

The organization of this work is as follows: Section \ref{sec: combinatorial}
considers the combinatorial minimization \eqref{opt: combinatorial}. 
Results regarding the continuous minimization \eqref{opt: continuous} 
are described in Section \ref{sec: continuous}. Finally, in Section \ref{sec: conclusion} we put forth and discuss a conjecture about an extension of the results of Section \ref{sec: continuous}. 

\section{Combinatorial minimization}\label{sec: combinatorial}
The main result of this section is stated in the forthcoming Theorem \ref{thm: main result}. For its statement, we need some definitions: 
\begin{definition}
Let $0<a<b$ be some integers and $I_{a,b}\equiv\{a,a+1,\ldots,b\}$. 
Then, $j \in I_{a,b}$ is a middle point of the set $I_{a,b}$ if and only if 
\begin{equation}
    \left\lfloor\frac{a+b}{2}\right\rfloor\leq j\leq\left\lceil\frac{a+b}{2}\right\rceil\,.
\end{equation}
\end{definition}
\begin{remark}\normalfont
Note that once $b-a$ is odd, there is a unique middle point of $I_{a,b}$. In addition, if $b-a$ is even, there are two middle points of $I_{a,b}$. 
\end{remark}

\begin{definition}\label{def: Gamma}
 Assume that $w \in (0,nx)$, $m\equiv\left\lfloor \frac{w}{x}\right\rfloor$ and $r\equiv w-mx$. 
 Then, for every positive integer $\delta \leq n+1$, let $\Gamma_\delta\equiv\Gamma_\delta(x,w)$ be the set of solutions $\boldsymbol{v}=(v_1,v_2,\ldots,v_n)$ with
\begin{equation}\label{eq: structure} v_i\equiv \begin{cases} 
      x & i\in\{\Delta_0,\Delta_0+\Delta_1,\ldots,\sum_{k=0}^{m-1}\Delta_k\} \\
      r & i=j \\
      0 & \text{o.w.} 
   \end{cases}
	\end{equation}
    such that:
\begin{itemize}
    \item $\Delta_0,\Delta_1,\ldots,\Delta_m$ are positive integers for which:
    \begin{enumerate}
        \item $ \sum_{k=0}^m\Delta_k=n+1$.
        \item $\max\left\{\Delta_k;0\leq k\leq m\right\}=\delta$.
        \item $t$ is some element of $\argmax\left\{\Delta_k;0\leq k\leq m\right\}$ and\newline $\left|\Delta_{k_1}-\Delta_{k_2}\right|\leq 1$, for every $k_1,k_2\in\{0,1,\ldots,m\}\setminus\{t\}$.
    \end{enumerate}
    
    \item  $j$ is a middle element of the set
	\begin{equation} \left\{\sum_{k=0}^{t-1}\Delta_k+1,\sum_{k=0}^{t-1}\Delta_k+2,\ldots,\sum_{k=0}^{t-1}\Delta_k+\delta-1\right\}\,.
	\end{equation}
\end{itemize}
\end{definition}

\begin{remark}\label{remark: m=0}
\normalfont When $w<x$, $m = 0$, and $\Gamma_\delta$ includes all solutions of the form \eqref{eq: structure} which are parameterized with 
$\Delta_0=\delta=n+1$, $t=0$, and $j$ which is a middle point 
of $\{1,2,\ldots,n\}$. 
\end{remark}

Note that the conditions $w \in (0, nx)$ and $x > 0$ imply that $n\geq1$. In particular, when $n=1$ we know that $w<x$, such that there is a unique optimal solution of \eqref{opt: combinatorial} which equals $\boldsymbol{v}=w$. Observe that,
as explained in Remark~\ref{remark: m=0},
it is the only element of $\Gamma_2$. 
Theorem~\ref{thm: main result} below refers to the case $n \geq 2$. 
\begin{definition} \label{def: h-function}
Let $a<b$ be two positive integers. In addition, assume that $\Delta_0(a,b),\Delta_1(a,b),\ldots,\Delta_{a-1}(a,b)$ are positive integers whose sum is $b$ and satisfy the condition:
 \begin{equation}\label{eq: equidistant}
 \max_{0\leq k_1,k_2\leq a-1}\left|\Delta_{k_1}(a,b)-\Delta_{k_2}(a,b)\right|\leq 1\,.
 \end{equation}
 Then, define
 \begin{align}\label{eq: h}
     h(a,b)&\equiv\frac{1}{2}\sum_{k=0}^{a-1}\Delta_k(a,b)\left[\Delta_k(a,b)-1\right]\\&=\frac{1}{2}\left\{\left(b\ \text{\normalfont mod}\ a\right)\left(\left\lfloor\frac{b}{a}\right\rfloor+1\right)\left\lfloor\frac{b}{a}\right\rfloor+\left[a-\left(b\ \text{\normalfont mod}\ a\right)\right]\left\lfloor\frac{b}{a}\right\rfloor\left(\left\lfloor\frac{b}{a}\right\rfloor-1\right)\right\}\,.\\&\nonumber\nonumber\end{align}
 \end{definition}
Now we are ready to state the main result:
\begin{theorem}\label{thm: main result}Assume that $w<nx$ and $n\geq2$. In addition, denote  $m\equiv\left\lfloor\frac{w}{x}\right\rfloor$ and $r\equiv w-mx$.
\begin{enumerate}
    \item If $m=0$, the set of optimal solutions of \eqref{opt: combinatorial} is $\normalfont\text{Conv}(\Gamma_{n+1})$ 
    (see also Remark \ref{remark: m=0}).
    
    \item If $m\geq1$, define two functions: 
    \begin{equation}\label{eq: A}
    A(\delta)\equiv A(\delta;r)\equiv x\left[h\left(m,n\!+\!1\!-\!\delta\right)+\frac{1}{2}\delta(\delta\!-\!1)\right]-r\left\lfloor\frac{\delta}{2}\right\rfloor\left\lceil\frac{\delta}{2}\right\rceil\,,
    \end{equation}  
    \begin{equation}
    \phi(\delta) \equiv 
    \left(2\delta+1\right)
    \left(x - \frac{r}{2}\right)
    -x
    \left(
    \left\lfloor \frac{n-\delta-1}{m}
    \right\rfloor
    +
    \left\lfloor \frac{n-\delta}{m}
    \right\rfloor
    \right) - \frac{r}{2}\,.
    \end{equation}
In addition, denote
\begin{align} \label{eq: delta1}
    \delta_1 &\equiv \min\{\delta \in \{1,3,5,\ldots\} \,;\, \phi(\delta) > 0\}\,, \\\nonumber
    \\
    \label{eq: delta2}
    \delta_2 &\equiv \min\{\delta \in \{2,4,6,\ldots\} \,;\, \phi(\delta) > 0\}\,.
\end{align}
and set
\begin{equation}\label{eq: delta*}
    \delta^*\equiv\begin{cases} 
      \delta_1 & \quad \text{\normalfont if }
      A(\delta_1) < A(\delta_2) \\
      \delta_2 & \quad \normalfont\text{otherwise}  
   \end{cases}
   \,.
\end{equation}
Then, $\delta^*$ is nondecreasing in $r$ and 
\begin{equation}
\label{eq: main result}
    \normalfont\text{Conv}\left(\Gamma_{\delta^*}\right)
    \subseteq\arg\min_{\boldsymbol{v}\in\Upsilon }f(\boldsymbol{v})\,.
\end{equation}
In particular, if $r=0$, then $\delta^*=1 + \left\lfloor \frac{n}{m+1} \right\rfloor = \left\lceil\frac{n+1}{m+1}\right\rceil$.
\end{enumerate}
\end{theorem}
\subsection{Implementation of Theorem \ref{thm: main result}}
When either $m=0$ or $r=0$, Theorem \ref{thm: main result} includes a precise construction of an optimal solution of the combinatorial minimization \eqref{opt: combinatorial} which is straightforward to implement. Otherwise, 
when $m\geq1$ and $r>0$, Theorem \ref{thm: main result} implies that the derivation of an optimal solution of \eqref{opt: combinatorial} boils down to the computations of $\delta_1$ and $\delta_2$ which are respectively given in 
\eqref{eq: delta1} and \eqref{eq: delta2}. In what follows, it is explained how to execute these computations.

To begin with, observe that $\delta_1$ and $\delta_2$ are
bounded by the solutions of the next two equations in $\delta$:
\begin{equation}
    \left(2\delta+1\right)\left(x-\frac{r}{2}\right)-x\left(\frac{n-\delta-1}{m}+\frac{n-\delta}{m}\right)-\frac{r}{2}=\pm2\,,
\end{equation}
i.e., if we apply the notations
\begin{equation}
    \delta_-\equiv\frac{\frac{r}{2}+x\frac{2n-1-m}{2m}-1}{x\left(1+\frac{1}{m}\right)-\frac{r}{2}}\  \ \ \text{and} \ \ \  \delta_+\equiv\frac{\frac{r}{2}+x\frac{2n-1-m}{2m}+1}{x\left(1+\frac{1}{m}\right)-\frac{r}{2}}\,,
\end{equation}
then we shall deduce that $\delta_-\leq\delta_1,\delta_2\leq\delta_+$.
This means that the search for $\delta_1$ (resp. $\delta_2$) may be limited to a set of integers 
\begin{align}
&\mathcal{D}_1\equiv\left\{\left\lceil\delta_-\right\rceil,\left\lceil\delta_-\right\rceil+1,\ldots,\left\lfloor\delta_+\right\rfloor\right\}\cap\left\{1,3,5,\ldots\right\}\,,\\&\nonumber\\&\nonumber \left(\text{resp. }\mathcal{D}_2\equiv\left\{\left\lceil\delta_-\right\rceil,\left\lceil\delta_-\right\rceil+1,\ldots,\left\lfloor\delta_+\right\rfloor\right\}\cap\left\{2,6,8,\ldots\right\}\right)\,.    
\end{align}
In particular, since $r\leq x$ and $m\geq1$, then 
\begin{equation}
     \delta_+-\delta_-=\frac{2}{x\left(1+\frac{1}{m}\right)-\frac{r}{2}}\leq\frac{4}{x}\,.
\end{equation} Therefore, we get an upper bound 
\begin{equation}
\max_{i=1,2}\left|\mathcal{D}_i\right|\leq
\frac{1}{2}\left\lceil\frac{4}{x}\right\rceil+1    
\end{equation}
which is fixed in all parameters except $x$. Furthermore, 
observe that
\begin{equation}
    \delta\mapsto\left(2\delta+1\right)\left(x-\frac{r}{2}\right)
    -x\left(\left\lfloor\frac{n-\delta-1}{m}\right\rfloor+\left\lfloor\frac{n-\delta}{m}\right\rfloor\right)
\end{equation}
is increasing in $\delta$. As a result,  for each $i=1,2$, $\delta_i$ may be derived via an application of a standard bisection method on the set $\mathcal{D}_i$.
Also, observe that the definition $\phi(\cdot)$ does not involve any complicated loops, and hence the application of the bisection method is relatively inexpensive. 

Lastly, once $\delta_1$ and $\delta_2$ are computed, \eqref{eq: delta*} dictates that we need to call $A(\cdot)$ twice in order to decide whether $\delta^*=\delta_1$ or $\delta^*=\delta_2$. Importantly, by definition (recall also the definition of $h(\cdot)$ in \eqref{eq: h}), the computation of $A(\cdot)$ involves only standard arithmetical operations. As a result, we may conclude that Theorem~\ref{thm: main result} yields an applicable numerical procedure to get an optimal solution of \eqref{opt: combinatorial}.
\subsection{Demonstration}
Before providing the proof of Theorem~\ref{thm: main result}, we would like to demonstrate some of its guidelines by an example. Specifically, consider the special case  $n=9$, $x = 1.1$, and $w = 2.4$. 
Each solution $\boldsymbol{v} = (v_1,\dots,v_{9})$ 
is characterized by a set $\{i_1,i_2,j\} \subseteq \{1,\dots,9\}$,
$i_1 < i_2$, such that $v_{i_1} = v_{i_2} = 1.1$, $v_{j} = 0.2$ and $v_i = 0$
otherwise.
\begin{figure}
    \centering
    \includegraphics{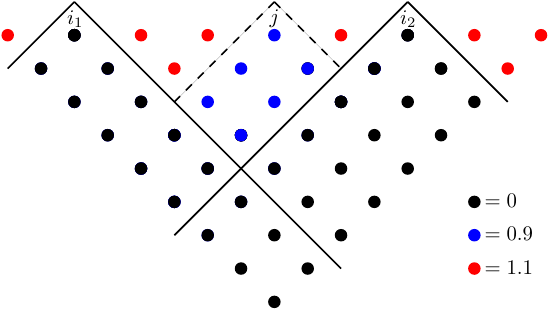}
    \caption{Representing $f(\cdot)$
    as the sum of the values of the individual balls when the parameters are $n=9$, $x=1.1$ and $w=2.4$.
    }
    \label{fig: triangle representation}
\end{figure}

Observe that by constructing a triangle
as in the example displayed in Figure~\ref{fig: triangle representation}, in which the value of the
$k$-th ball in row $l$ is given by $x - \sum_{i = k}^{l + k} v_i$,
for each choice of $\{i_1,i_2,j\}$
the value of 
$f(\boldsymbol{v})$ can be found by
summing over the values of the individual balls in
such a triangle. Notably, for each choice of $i_1,i_2$,
there are at most three sub-triangles whose mass is unequal
to $0$, and within one of these triangles there
is a rectangle within which the
balls attain mass $0.9$. Therefore, 
with $\Delta_1 = i_1, \Delta_2 = i_2-i_1$
and $\Delta_3 = 10-i_2$, 
for $i_1 < j < i_2$, we may write 
\begin{align*}
    f(\boldsymbol{v}) = 1.1 \cdot
    \left[\frac{\Delta_1
    (\Delta_1\!-\!1)}{2}\!+\!\frac{\Delta_2
    (\Delta_2\!-\!1)}{2}\!+\!\frac{\Delta_3
    (\Delta_3\!-\!1)}{2}\right] - 0.2 \cdot (j\!-\!i_1)(i_2\!-\!j).
\end{align*}
Now, given $i_1$ and $i_2$, observe
that the number of balls within
the rectangular area is maximal if $j$ is located
at a middle point of the hypotenuse of the largest triangle 
(note that there are two middle points 
when the hypotenuse consist of an even number of balls).
Thus, $i_1,i_2$ should minimize
\begin{align*}
    f(\boldsymbol{v}) = 1.1 \cdot
    \sum_{i = 1}^3 \frac{\Delta_i
    (\Delta_i\!-\!1)}{2}  - 0.2 \cdot
    \left\lfloor\frac{\Delta_{t}}{2}\right\rfloor\left\lceil\frac{\Delta_{t}}{2}\right\rceil,
\end{align*}
with $\Delta_t = \max\{\Delta_1,\Delta_2,\Delta_3\}$.
\subsection{Proof}
Before providing the proof, we need to discuss some preliminaries. Namely, notice that in Definition~\ref{def: h-function}
it is always possible to choose $\{\Delta_k(a,b)\}_{k = 0}^{a-1}$
with $\Delta_0(a,b)=\left\lfloor\frac{b}{a}\right\rfloor$. Thus, it is possible to choose $\Delta_0(a,b+1),\Delta_1(a,b+1),\ldots,\Delta_{a-1}(a,b+1)$ such that for every $0\leq k\leq a-1$:
\begin{equation} \Delta_k(a,b+1)\equiv \begin{cases} 
      \Delta_{0}(a,b)+1 & k=0\,, \\
      \Delta_k(a,b) & \text{o.w.} 
   \end{cases}
	\end{equation}
and hence
	\begin{align}\label{eq: h-dif}
	2\left[h(a,b+1)-h(a,b)\right]&=\left[\Delta_0(a,b)+1\right]\Delta_0(a,b)-\Delta_0(a,b)\left[\Delta_0(a,b)-1\right]\nonumber\\&=2\Delta_0(a,b)=2\left\lfloor\frac{b}{a}\right\rfloor\,.
	\end{align}
    Another useful observation is that for any positive integer $\delta$,
	\begin{equation}
	4\left\lfloor\frac{\delta}{2}\right\rfloor\left\lceil\frac{\delta}{2}\right\rceil=\begin{cases} 
      \delta^2 & \text{$\delta$ is even} \\
      \delta^2-1 & \text{$\delta$ is odd} 
   \end{cases}
	\end{equation}
	and hence 
	\begin{align}\label{eq: eta-value}
	 \eta(\delta)&\equiv\left\lfloor\frac{\delta+1}{2}\right\rfloor\left\lceil\frac{\delta+1}{2}\right\rceil-\left\lfloor\frac{\delta}{2}\right\rfloor\left\lceil\frac{\delta}{2}\right\rceil=\frac{1}{2}\begin{cases} 
      \delta & \text{$\delta$ is even} \\
      \delta+1 & \text{$\delta$ is odd} 
   \end{cases}\,.
	\end{align}
Now, we are ready to provide the proof of Theorem \ref{thm: main result}. 
\newline\newline
\begin{proof}
For every $\{i_1,i_2,\ldots,i_m,j\}\subseteq\{1,2,\ldots,n\}$ define 	\begin{equation} v_i(i_1,i_2,\ldots,i_m,j)\equiv \begin{cases} 
      x & i\in\{i_1,i_2,\ldots,i_m\} \\
      r & i=j \\
      0 & \text{o.w.} 
   \end{cases}\ \ , \ \ 1\leq i\leq n
	\end{equation}    
	and the resulting  solution	
	\begin{equation}
	    \boldsymbol{v}(i_1,i_2,\ldots,i_m,j)=\left(v_i(i_1,i_2,\ldots,i_m,j)\right)_{1\leq i\leq n}\,.
	\end{equation}
	Without loss of generality, assume that $i_1<i_2<\ldots<i_m$ (otherwise sort these indices and rename them properly). In addition, let $(i_0,i_{m+1})\equiv(0,n+1)$ and for each $0\leq k\leq m$ define $\Delta_k\equiv i_{k+1}-i_k$. Furthermore, denote $t\equiv \max\{k\geq0;i_k\leq j\}$ and observe that
	the objective function is given by:
	\begin{align}\label{eq: objective}
     g(i_1,i_2,\ldots i_m,j)&\equiv f\left[\boldsymbol{v}(i_1,i_2,\ldots,i_m,j)\right]\\&=x\sum_{k=0}^{m}\psi\left(\Delta_k\right)-r(j-i_t)(i_{t+1}-j)\,,\nonumber
	\end{align}
	where $\psi(y)\equiv\frac{y(y-1)}{2}$, $y\geq1$.
	Now, given fixed  $i_1,i_2,\ldots,i_m$, then \eqref{eq: objective} reveals that the optimal $j$  must be a middle point of $\{i_{t}+1,i_{t}+2,\ldots,i_{t+1}-1\}$. Notably, once there are two middle points, a symmetry argument implies that they yield the same objective value and hence they are equivalent. This observation completes the proof for the first theorem's statement due to the convexity of $f(\cdot)$.
	
	In order to complete the proof we need to assume that $m\geq1$ and it is left to minimize
	\begin{align}
	 \widetilde{g}&(\Delta_0,\Delta_1,\ldots,\Delta_{m})\equiv x\sum_{k=0}^{m}\psi(\Delta_k)-r\left\lfloor\frac{\Delta_{t}}{2}\right\rfloor\left\lceil\frac{\Delta_{t}}{2}\right\rceil\,.
	\end{align}
	In particular, $\psi(\cdot)$ is strictly convex and increasing on $[1,\infty)$. Consequently, the optimal $\{\Delta_k;0\leq k\leq m\}$ is characterized by the next conditions:
	\begin{enumerate}
	    \item $\Delta_0+\Delta_1+\ldots+\Delta_m=n+1$.
	    
	    \item $\Delta_t=\max\left\{\Delta_k;0\leq k\leq m\right\}$ and\newline $\left|\Delta_{k_1}-\Delta_{k_2}\right|\leq 1$, for every $k_1,k_2\in\{0,1,\ldots,m\}\setminus\{t\}$.
	    
	\item $\Delta_t$ is a minimizer of the function $A(\delta)$ (in $\delta)$.
    \end{enumerate}
	Especially, \eqref{eq: h-dif} and \eqref{eq: eta-value} imply that
	\begin{align}
    {\rm d}_1A(\delta)&\equiv 
    A(\delta+1)-A(\delta)
    \\
    &\nonumber=x\left(\delta-\left\lfloor\frac{n-\delta}{m}\right\rfloor\right)-r\eta(\delta)\\&=\delta\left(x-\frac{r}{2}\right)-x\left\lfloor
    \frac{n-\delta}{m}\right\rfloor-\frac{r\textbf{1}_{\{1,3,5,\ldots\}}(\delta)}{2}\,.\nonumber
	\end{align}
Consequently, deduce that
     \begin{align}
      {\rm d}_2A(\delta)&\equiv 
    A(\delta+2)-A(\delta)
  ={\rm d}_1A(\delta+1) + {\rm d}_1A(\delta)
    = \phi(\delta)\,.
    \end{align}
Since ${\rm d_2}A(\cdot)$ is increasing in $\delta$, the minimizer of $A(\cdot)$ must be $\delta^*$ as given in the statement of the theorem, and hence \eqref{eq: main result} stems from the convexity of $f(\cdot)$. 
 
When $r=0$,
 ${\rm d_1}A(\cdot)$ is also increasing in $\delta$,
 such that $\delta^*$ is equivalently characterized by the minimum
 $\delta > 0$ for which ${\rm d_1}A(\delta) > 0$, which yields $\Delta_t$ for this case. It is left to prove that $\delta^*$ is nondecreasing in $r$. To this end, take $0\leq r_1<r_2<x$ and notice that for every $\delta$, the difference
	\begin{equation}
	    A(\delta;r_2)-A(\delta;r_1)=-(r_2-r_1)\left\lfloor\frac{\delta}{2}\right\rfloor\left\lceil\frac{\delta}{2}\right\rceil
	\end{equation}
	is decreasing in $\delta$. Thus, Theorem 10.6 in \cite{Sundaram1996} implies the required result.  
\end{proof}
\section{Continuous minimization}\label{sec: continuous}
\subsection{When $r=0$}
The forthcoming Theorem \ref{thm: continuous} reveals that once $r=0$, the class of optimal solutions of \eqref{opt: combinatorial} which is stated in Theorem \ref{thm: main result} is actually optimal in the broader (compact) minimization domain $\Lambda$ defined in \eqref{eq: lambda}. Before providing its statement, we present the following useful notations:
\begin{equation}
\tau_u \equiv \tau_u(n,m) \equiv \left\lceil \frac{n+1}{m+1} \right\rceil\ \ , \ \ \tau_{l} \equiv \tau_l(n,m) \equiv \left\lfloor \frac{n+1}{m+1} \right\rfloor\,.    
\end{equation} 
\begin{theorem}\label{thm: continuous}
    If $r=0$, i.e., $w=mx$ for some positive integer $m$, then,
\begin{equation}\normalfont\text{Conv}\left(\Gamma_{\tau_u}\right)\subseteq\arg\min_{\boldsymbol{v}\in\Lambda}f(\boldsymbol{v})
\end{equation}
and
\begin{equation}
    \min_{\boldsymbol{v}\in\Lambda}f(\boldsymbol{v})=\left(\tau_u-1\right)\left(x(n+1)-\frac{1}{2}(w+x)\tau_u\right)
\end{equation}
\end{theorem}

\begin{remark}\label{remark: efficient}
    \normalfont Note that construction of some vector in $\Gamma_{\tau_u}$ is not computationally demanding. Therefore, when $r=0$, Theorem~\ref{thm: continuous} yields a natural numerical procedure which is easy to implement in order to solve the continuous minimization \eqref{opt: continuous}.
\end{remark}

\begin{figure}[t]
    \centering
    \includegraphics{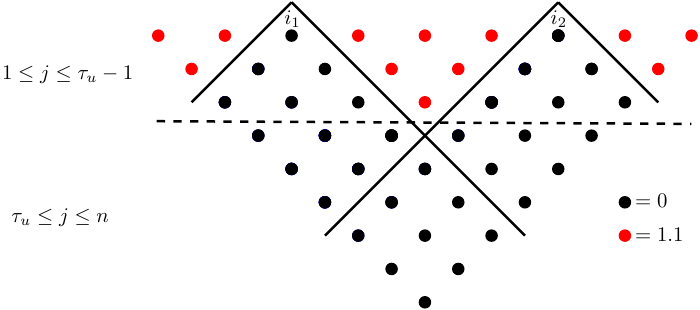}
    \caption{Representing $f(\cdot)$
    as the sum of the values of the individual balls when the parameters are $n=9$, $x=1.1$ and $w=2.2$. For each $1 \leq j \leq n$, the function $f_j(\cdot)$ is
    the sum of the balls in row $j$. The optimal solution spreads the mass in
    almost equal gaps, i.e., a minimum gap of $\tau_l = 3$, and a maximum gap of
    $\tau_u = 4$. For each $j > \tau_u$, the balls in row $j$ have no mass, whereas for each row $1 \leq j \leq \tau_l$, there is no overlap between the balls affected by the different vector points with mass.}
    \label{fig: triangle representation proof}
\end{figure}

\begin{proof}
Let
$\boldsymbol{v}^* $ be an element of $\Gamma_{\tau_u}$. In addition, for any $1\leq j\leq n$ and $\boldsymbol{v}\in\Lambda$ define
 \begin{equation}
    L_j\equiv\left\{(k,\ell)\in\mathbb{N}^2\ ;\ k\leq\ell\leq n\ , \ \ell-k+1=j\right\}\,,
\end{equation}
and
 \begin{align} \label{def: f_j}
        f_j(\boldsymbol{v}) \equiv \sum_{(k,\ell)\in L_j} 
        \left(x - \sum_{i=k}^{\ell} v_i\right)^+\,,
    \end{align}
such that
\begin{align} \label{eq: totalsum}
    f(\boldsymbol{v}) = \sum_{j = 1}^n f_j(\boldsymbol{v})\,.
\end{align}
In what follows, we will prove that 
$f_j(\boldsymbol{v}) \geq f_j(\boldsymbol{v}^*)$ for any 
$1\leq j\leq n$ and $\boldsymbol{v}\in\Lambda$, thereby
providing an explicit expression
for $f_j(\boldsymbol{v^*})$. 
This will end the proof since $f(\cdot)$ is convex,
and the minimum objective function is found through
\eqref{eq: totalsum} with input $\boldsymbol{v}^*$. 

To begin with, as illustrated in in Figure~\ref{fig: triangle representation proof},
observe that
\begin{equation}
    \max\left\{i_2-i_1;1\leq i_1<i_2\leq n\ \ \text{s.t.}\ v_{i_1}^*=v_{i_2}^*=x, v_{i_1+1}^*=\ldots=v_{i_2-1}^* = 0\right\}
\end{equation}
is bounded from above by $\tau_u$. Consequently, $f_j(\boldsymbol{v}^*)=0$ for every $j\geq\tau_u$, i.e., it is left to consider $1 \leq j \leq \tau_u-1$. To this end, consider some $1\leq j\leq \tau_u-1$ and recall that 
\begin{equation}\label{eq: min gap}
    \min\left\{i_2-i_1;1\leq i_1<i_2\leq n\ \ \text{s.t.}\ v_{i_1}^*=v_{i_2}^*=x\right\}=\tau_l \geq \tau_u-1\,.
\end{equation}
Thus, since $\sum_{i=k}^\ell v_i^*\leq x$ for every $(k,\ell)\in L_j$, deduce that
\begin{equation}
 \left|\left\{(k,\ell)\in L_j\ ;\ \sum_{i=k}^\ell v_i^*=x\right\}\right|=jm
\end{equation}
and hence
\begin{align}\label{eq: f_j equidistant}
    f_{j}(\boldsymbol{v}^*)&= 
    \sum_{(k,\ell)\in L_j} 
        \left(x - \sum_{i=k}^{\ell} v_i^*\right)\\
        &\nonumber=(n\!+\!1\!-\!j)x-\sum_{(k,\ell)\in L_j}\sum_{i=k}^{\ell}v_i^*\\
        &\nonumber=(n\!+\!1\!-\!j)x-jmx\\& = (n\!+\!1\!-\!j)x-jw\,.\nonumber
\end{align}
Now, consider some $\boldsymbol{v}\in\Lambda(m)\cap[0,x]^n$ and observe that
\begin{align}
    f_j(\boldsymbol{v})&
    \geq \sum_{(k,\ell)\in L_j}\left(x-\sum_{i=k}^\ell v_i\right)
    =(n+1-j)x-\sum_{(k,\ell)\in L_j}\sum_{i=k}^{\ell}v_i\,.
\end{align}
Note that for each $1\leq i\leq n$ we have 
\begin{equation}
    \left|\left\{(k,\ell)\in L_j;k\leq i\leq\ell\right\}\right|\leq j
\end{equation}
and hence
\begin{equation}
    \sum_{(k,\ell)\in L_j}\sum_{i=k}^{\ell}v_i\leq jw\,,
\end{equation}
from which the result follows. 
\end{proof}

\subsection{When $0<r<x$}
It is natural to check whether Theorem~\ref{thm: continuous} can be extended to the case $0<r<x$. Example~\ref{example: suboptimal} 
below demonstrates that such a generalization is not straightforward. 

\begin{example}\label{example: suboptimal}\normalfont
Consider the parameterization $n = 7$, $x = 1$, and $w = 2.2$. 
According to Theorem~\ref{thm: main result},
one optimal solution of \eqref{opt: combinatorial} is given by $\boldsymbol{v}_1=(0,1,0,0,1,0.2,0)$, with $f(\boldsymbol{v}_1) = 6.6$. 
However, $\boldsymbol{v}_1$ is a sub-optimal solution of \eqref{opt: continuous}, as the vector 
$\boldsymbol{v}_2=(0,0.2,0.8,0.2,0.8,0.2,0)$ yields $f(\boldsymbol{v}_2) = 6.4$. 
\end{example}
Notably, in Example \ref{example: suboptimal}, 
$\boldsymbol{v}_2$ can be
written as a convex combination of
$\boldsymbol{v} = (0,0,1,0,1,0,0) \in \Gamma_{\tau_u(n,m)}$
and $\tilde{\boldsymbol{v}} = (0,1,0,1,0,1,0) \in 
\Gamma_{\tau_u(n,m+1)}$,
which are optimal solutions
when $w = mx$ and $w = (m+1)x$, respectively. 
Therefore, we still get a reduction of the continuous minimization \eqref{opt: continuous} to the combinatorial one \eqref{opt: combinatorial}.
The following Theorem \ref{thm: conjecture} expands this observation to a wider family of parameterizations.

\begin{theorem}\label{thm: conjecture}
    Denote $y\equiv y(x,w)\equiv x-r$,
    $\tau_{u,1} \equiv \tau_u(n,m)$, 
    $\tau_{l,1} \equiv \tau_l(n,m)$,
    $\tau_{u,2} \equiv \tau_u(n,m+1)$,
    and $\tau_{l,2} \equiv \tau_l(n,m+1)$,
    and consider any two vectors 
    $\boldsymbol{v}_y\in\Gamma_{\tau_{u,1}}\left(y,my\right)$
    and $\boldsymbol{v}_r\in\Gamma_{\tau_{u,2}}\left(r,\left(m+1\right)r\right)$. 
    If $\tau_{u,1} = \tau_{u,2}$ or
    $\tau_{l,1} = \tau_{l,2}$, then
    \begin{equation}
\boldsymbol{v}^*\equiv\boldsymbol{v}^*\left(\boldsymbol{v}_y,\boldsymbol{v}_r\right)\equiv\boldsymbol{v}_y+\boldsymbol{v}_r
    \end{equation}
     is such that 
     \begin{equation}
     f(\boldsymbol{v}^*) = \min_{\boldsymbol{v}\in\Lambda}f(\boldsymbol{v})=\left(\tau_{u,1}-1\right)\left(x(n+1)-\frac{1}{2}(w+x)\tau_{u,1}\right)  \end{equation}
\end{theorem}

\begin{remark}
    \normalfont Prominently, the objective values in Theorems~\ref{thm: continuous}
    and~\ref{thm: conjecture} are given by the same formula. Considering a
    representation as in Figure~\ref{fig: triangle representation proof}, 
    an optimal vector from Theorem~\ref{thm: continuous} 
    allocates the elements with mass such
    that, before row $\tau_u$, no elements with mass affect the same balls,
    and from row $\tau_u$ onward, all balls are zero. A similar
    idea is true for an optimal vector in Theorem~\ref{thm: conjecture}:
    No two elements whose mass sums to $x$ affect the same balls 
    before $\tau_u$, and from $\tau_u$ onward, all balls have zero mass.
\end{remark}

\begin{example}
    \normalfont Theorem~\ref{thm: conjecture} is \textit{not} an extension of Theorem~\ref{thm: continuous}. For example,
    if $n = 7, x = 1$, and $w = 2.2$, 
    then \textit{only} Theorem \ref{thm: conjecture} can be applied. To see this, observe that
    $m = \left\lfloor \frac{2.2}{1} \right\rfloor = 2$
    and $r = 2.2 - 2 \cdot 1 > 0$, i.e., it is impossible to apply Theorem \ref{thm: continuous}. On the other hand, notice that
    \begin{equation*}
        \tau_{l}(n,m) = \left\lfloor \frac{8}{3} \right\rfloor 
        = 2 = \left\lfloor \frac{8}{4} \right\rfloor = \tau_{l}(n,m+1).    \end{equation*}
        and hence Theorem \ref{thm: conjecture} may be applied.
\end{example}

\begin{example}
\normalfont    
Theorem~\ref{thm: continuous} is \textit{not} an extension of Theorem~\ref{thm: conjecture}. For example,  if $n = 8, x = 1$, and $w = 3$, then \textit{only} Theorem \ref{thm: continuous} can be applied. To see this, observe that
    $m = \left\lfloor \frac{3}{1} \right\rfloor = 3$
    and $r = 3 - 3 \cdot 1 = 0$, i.e., Theorem \ref{thm: continuous} may be applied. On the other hand, notice that
    \begin{align}
        \tau_{l}(n,m) &= \left\lfloor \frac{9}{4} \right\rfloor 
        \neq \left\lfloor \frac{9}{5} \right\rfloor = \tau_{l}(n,m+1) 
        \quad , \quad
        \tau_{u}(n,m) = \left\lceil \frac{9}{4} \right\rceil 
        \neq \left\lceil \frac{9}{5} \right\rceil = \tau_{u}(n,m+1)
        \end{align}
        and hence it is impossible to apply Theorem \ref{thm: conjecture}.
\end{example}

\begin{proof}
    Define
    $L_j$ and $f_j(\boldsymbol{v})$
    as in the proof of Theorem~\ref{thm: continuous}.
    Then, observe that for each $z \in \{y,r\}$
    \begin{equation}
    \max\left\{i_2-i_1;1\leq i_1<i_2\leq n\ \ \text{s.t.}\ v_{z,i_1}=v_{z,i_2}=z, v_{z,i_1+1}=\ldots=v_{z,i_2-1}=0 \right\}
\end{equation}
is bounded from above by $\tau_{u,1}$. Therefore, the fact that $y+r=x$ implies that $f_j(\boldsymbol{v}^*) = 0$ for every $j \geq \tau_{u,1}$. Moreover, observe that 
    for each $z \in \{y,r\}$,
    \begin{equation}\label{eq: min gap 2}
    \min\left\{i_2-i_1;1\leq i_1<i_2\leq n\ \ \text{s.t.}\ v_{z,i_1}=v_{z,i_2}=z\right\} \geq\tau_{l,2}\,.
\end{equation}
In case that $\tau_{u,1} = \tau_{u,2}$, 
    then $\tau_{l,2} \geq \tau_{u,2} - 1 = \tau_{u,1}-1$.
    In addition, if $\tau_{l,1} = \tau_{l,2}$, then we get $\tau_{l,2} = \tau_{l,1} \geq \tau_{u,1}-1$.
    Thus, for $1 \leq j \leq \tau_{u,1}-1$ we shall apply 
    the same arguments as in the proof of Theorem~\ref{thm: continuous},
    together with the fact that $y+r=x$. Altogether, this implies that
    for any $\boldsymbol{v} \in \Lambda \cap [0,x]^n$
\begin{align}
    f_{j}(\boldsymbol{v}^*)
    &= 
    \sum_{(k,\ell)\in L_j} 
        \left(x - \sum_{i=k}^{\ell} v_i^*\right)\\
        &\nonumber=(n\!+\!1\!-\!j)x-\sum_{(k,\ell)\in L_j}\sum_{i=k}^{\ell}v_{y,i}
        -\sum_{(k,\ell)\in L_j}\sum_{i=k}^{\ell}v_{r,i}
        \\
        &\nonumber=(n\!+\!1\!-\!j)x-jmy - j(m\!+\!1)r\\
        &\nonumber=(n\!+\!1\!-\!j)x-jw
        \leq f_j(\boldsymbol{v})\,.
\end{align}
Finally, the objective value of $\boldsymbol{v}^*$ is found through \eqref{eq: totalsum}.
\end{proof}

\section{Conjecture}\label{sec: conclusion}
It is natural to look for a unifying generalization of Theorem \ref{thm: continuous} and \ref{thm: conjecture}. The following Conjecture \ref{conjecture} has been verified in an enormous number of numerical examples using Mathematica. 

\begin{conjecture}\label{conjecture}
    Denote $y\equiv y(x,w)\equiv x-r$,
    $\tau_1 \equiv \tau_u(n,m)$, $\tau_2 \equiv \tau_u(n,m+1)$,
    and
    consider any two vectors $\boldsymbol{v}_r\in\Gamma_{\tau_1}\left(r,\left(m+1\right)r\right)$ and $\boldsymbol{v}_y\in\Gamma_{\tau_2}\left(y,my\right)$. 
    Then, 
    \begin{equation}
\boldsymbol{v}^*\equiv\boldsymbol{v}^*\left(\boldsymbol{v}_y,\boldsymbol{v}_r\right)\equiv\boldsymbol{v}_y+\boldsymbol{v}_r
    \end{equation}
     is an optimal solution of \eqref{opt: continuous}.
\end{conjecture}

\begin{remark}
    \normalfont Observe that when $r=0$, Conjecture \ref{conjecture} becomes Theorem \ref{thm: continuous}.
    Moreover, note that the case $m \geq \lfloor \frac{n}{2} \rfloor$ is trivial,
    as $f_j(\boldsymbol{v}^*) = 0$ for any $j \geq 0$, and $f_1(\boldsymbol{v})$
    has the same value for each $\boldsymbol{v} \in \Lambda \cap [0,x]^n$.
    Therefore, Theorem \ref{thm: conjecture} implies that it is left to prove Conjecture \ref{conjecture} under the assumptions $0<r<x$,
    $m < \lfloor \frac{n}{2} \rfloor$, and $\tau_u(n,m) \neq \tau_u(n,m+1)$.
\end{remark}

\begin{remark}
    \normalfont Note that if Conjecture \ref{conjecture} is true, the convexity of $f(\cdot)$ and the convexity of $\Lambda$ imply that the resulting convex hull is a subset of optimal solutions. 
\end{remark}

\begin{remark}
    \normalfont Observe that the construction of $\boldsymbol{v}^*$ is based on optimal solutions of the combinatorial minimization \eqref{opt: combinatorial}. Therefore, we conjecture that the reduction of the continuous minimization \eqref{opt: continuous} to the combinatorial minimization \eqref{opt: combinatorial} is valid for any parameterization. 
\end{remark}
Prominently, the incomplete part of this research is the proof of Conjecture \ref{conjecture}. We anticipate that its proof should follow from (at least) one of the following directions: (1)~Phrasing the problem as a dynamic program which can be solved by an application of Theorem \ref{thm: continuous}, or (2)~applying the theory about Lagrange multipliers in the context of sub-differential calculus (see, e.g., \cite{Rockafellar1997}). Our experience teaches that while both of these approaches look reasonable, they are not straightforward to apply. Therefore, we consider providing the proof of Conjecture~\ref{conjecture} as an open problem of different magnitude
that goes beyond the scope of the current research. Hopefully, it will be solved in the future.
\newline\newline
\textbf{Acknowledgement:} We would like to thank Danny Segev and Leen Stougie for their comments on earlier versions of the current manuscript. 
\newline\newline
\textbf{Funding:} This research was supported by the European Union’s Horizon 2020 research and innovation
programme [Marie Skłodowska-Curie Grant Agreement 945045] and the NWO Gravitation project
NETWORKS [Grant 024.002.003].

\end{document}